\documentclass[reqno]{amsart}
\usepackage{amsmath}
\usepackage{amssymb}
\usepackage{amsfonts}

\newtheorem{theorem}{Theorem}
\theoremstyle{plain}

\newtheorem{claim}{Claim}

\newtheorem{corollary}{Corollary}

\newtheorem{definition}{Definition}
\newtheorem{example}{Example}

\newtheorem{proposition}{Proposition}
\newtheorem{remark}{Remark}

\numberwithin{equation}{section}

\begin{document}
\title[Uniform disconnectedness and Quasi-Assouad Dimension]{Uniform
disconnectedness and Quasi-Assouad Dimension}
\author{Fan L\"u}
\address{School of Mathematics and Statistics, Huazhong University of
Science and Technology, 430074, Wuhan, P. R. China}
\email{lvfan1123@163.com}
\author{Li-Feng Xi}
\address{Institute of Mathematics, Zhejiang Wanli University, Ningbo,
Zhejiang, 315100, P.~R. China}
\email{xilifengningbo@yahoo.com}
\thanks{Li-Feng Xi is the corresponding author. This work is supported by
National Natural Science Foundation of China (Nos. 11371329, 11071224) and
NSF of Zhejiang Province (Nos. LR13A1010001, LY12F02011)}

\begin{abstract}
The uniform disconnectedness is an important invariant property under
bi-Lipschitz mapping, and the Assouad dimension $\dim _{A}X<1$ implies the
uniform disconnectedness of $X$. According to quasi-Lipschitz mapping, we
introduce the quasi-Assouad dimension $\dim _{qA}$ such that $\dim _{qA}X<1$
implies its quasi uniform disconnectedness. We obtain $\overline{\dim }%
_{B}X\leq \dim _{qA}X\leq \dim _{A}X$ and compute the quasi-Assouad
dimension of Moran set.
\end{abstract}

\maketitle

\section{Introduction}

A subset $E$ of metric space $(X,\mathrm{d})$ is said to be \textbf{%
uniformly disconnected \cite{David Semmes}}, if there is a constant $0<c<1$
such that for any $x\in E\ $and any $0<r<r^{\ast }$ for some $r^{\ast },$
there exists a set $E_{x,r}\subset E$ satisfying
\begin{equation}
E\cap B(x,cr)\subset E_{x,r}\subset B(x,r)\text{ and }\mathrm{dist}%
(E_{x,r},E\backslash E_{x,r})\geq cr,  \label{d:UD}
\end{equation}%
where $B(x,r)$ is the closed ball with center $x$ and radius $r,$ and $\,%
\mathrm{dist}(\cdot ,\cdot )$ denotes the least distance between sets. This
uniform disconnectedness is an invariant property under any bi-Lipschitz
mapping. Here $f:(X,\mathrm{d}_{X})\rightarrow (Y,\mathrm{d}_{Y})$ is
bi-Lipschitz$,$ if there exists a constant $L>0$ such that for all $%
x_{1},x_{2}\in X,$
\begin{equation}
L^{-1}\mathrm{d}_{X}(x_{1},x_{2})\leq \mathrm{d}_{Y}(f(x_{1}),f(x_{2}))\leq L%
\mathrm{d}_{X}(x_{1},x_{2}).
\end{equation}

The uniform disconnectedness plays an important role in fractal geometry.
David and Semmes \cite{David Semmes} obtained the uniformization theorem on
quasisymmetric equivalence: if a compact metric spaces is uniformly
disconnected, uniformly perfect and doubling, then they are
quasisymmetrically equivalent to the Cantor ternary set $C$. Mattila and
Saaranen \cite{MS} proved that suppose $E$ and $F$ are Ahlfors-David $s$%
-regular and $t$-regular respectively with $s<t,$ if $E$ is uniformly
disconnected, then $E$ can be bi-Lipschitz embedded into $F.$ Based on the
method of \cite{MS}, Wang and Xi \cite{Wang Xi N} discussed the
quasi-Lipschitz equivalence between uniformly disconnected Ahlfors-David
regular sets.

\medskip

An interesting property (e.g. see Proposition 5.1.7 of \cite{MT}) is that
\begin{equation}
\dim _{A}X<1\implies X\text{ is uniformly disconnected.}
\end{equation}%
We can recall \textbf{Assouad dimension} $\dim _{A}$ as follows. We say $(X,%
\mathrm{d})$ is \textit{doubling} if there exists an integer $N>0$ such that
each closed ball in $X$ can be covered by $N$ closed balls of half the
radius. Repeated applying the doubling property, it gives that there exist
constants $b,c>0$ and $\alpha >0$ such that for all $r$ and $R$ with $%
0<r<R<b $, every closed ball $B(x,R)$ can be covered by $c(\frac{R}{r}%
)^{\alpha }$ balls of radius $r$. Let $N_{r,R}(X)$ denote the smallest
number of balls with radii $r$ needed to cover any ball with radius $R$. The
Assouad dimension of $X$, denoted by $\dim _{A}X$, is defined as
\begin{equation*}
\dim _{A}X=\inf \{\alpha \geq 0\text{ }|\text{ }\exists \text{ }b,c>0\text{
s.t. }N_{r,R}(X)\leq c(\frac{R}{r})^{\alpha }\text{ }\forall \text{ }%
0<r<R<b\},
\end{equation*}%
which was introduced by Assouad in the late 1970s \cite%
{Assouad1,Assouad2,Assouad3}. Now it plays a prominent role in the study of
\textbf{quasiconformal mappings} and \textbf{embeddability problems}, and we
refer the readers to the textbook \cite{Heinonen} and the survey paper \cite%
{Luukkainen} for more details. Olsen \cite{Olsen} obtained the Assouad
dimensions for a class of fractals with some flexible graph-directed
construction, Mackay \cite{Mackay} and Fraser \cite{Fraser} calculated the
Assouad dimensions of some classes of self-affine fractals. It is well known
that $\dim _{H}X\leq \overline{\dim }_{B}X\leq \dim _{A}X,$ where $\dim
_{H}(\cdot )$ and $\overline{\dim }_{B}(\cdot )$ are Hausdorff and upper box
dimensions respectively. For example. if $E$ is Ahlfors-David $s$-regular
\cite{David Semmes}$,$ then $\dim _{A}E=\dim _{H}E=s,$ furthermore, if $s<1$
then $E$ is uniformly disconnected.

We say that a bijection $f:X\rightarrow Y$ is a \textbf{quasi-Lipschitz }%
mapping, if for all $x_{1},x_{2}\in X,$%
\begin{equation}
\frac{\log \text{d}_{Y}(f(x_{1}),f(x_{2}))}{\log \text{d}_{X}(x_{1},x_{2})}%
\rightarrow 1\text{ uniformly as }\mathrm{d}_{X}(x_{1},x_{2})\rightarrow 0.
\label{qqq}
\end{equation}%
We say $X$ and $Y$ are quasi-Lipschitz equivalent, if the above
quasi-Lipschitz mapping exists. We can introduce quasi-H\"{o}lder
equivalence and mapping, if $1$ in (\ref{qqq}) is replaced by $\dim
_{H}X/\dim _{H}Y.$

Inspired by \cite{MS}, Wang and Xi \cite{Wang Xi S} introduced the quasi
Ahlfors-David regularity and \textbf{quasi} \textbf{uniform disconnectedness}
and proved that $E$ is quasi-H\"{o}lder equivalent to the Cantor ternary set
$C$ if and only if $E$ is quasi Ahlfors-David regular and quasi uniformly
disconnected. As a consequence, they obtained that if quasi Ahlfors-David
regular sets $E$ and $F$ are quasi uniformly disconnected, then $E$ and $F$
are quasi-Lipschitz equivalent if and only if $\dim _{H}E=\dim _{H}F$.

\begin{definition}
We say that a compact subset $E$ of metric space $X$ is quasi uniformly
disconnected, if there is a function $\psi :\mathbb{R}^{+}\rightarrow
\mathbb{R}^{+}$ with $\psi (r)<r$ for $r>0$ and $\lim_{r\rightarrow 0}\frac{%
\log \psi (r)}{\log r}=1$ such that for any $x\in E\ $and any $0<r<r^{\ast }$
$(r^{\ast }$ is a constant$),$ there exists a set $E_{x,r}\subset E$
satisfying%
\begin{equation}
E\cap B(x,\psi (r))\subset E_{x,r}\subset B(x,r)\text{ and }\mathrm{dist}%
(E_{x,r},E\backslash E_{x,r})\geq \psi (r).  \label{d:QUD}
\end{equation}
\end{definition}

Note that the quasi uniform disconnectedness is an invariant property under
any quasi-H\"{o}lder(Lipschitz) mapping.

\begin{example}
Given $0<\alpha <1$. Let $a_{k}=\prod_{i=1}^{k}(1-(i+1)^{-\alpha })$ for all
$k\geq 1$. We can check that the countable set $E=\{0,1,a_{1},a_{2},\cdots
\} $ is quasi uniformly disconnected but not uniformly disconnected. We give
the detail of this example in Section 2.
\end{example}

\subsection{Quasi-Assouad dimension}

$\ $

The motivation of this manuscript is to introduce a notion named
quasi-Assouad dimension $\dim _{qA}X$ satisfying that
\begin{equation*}
\dim _{qA}X<1\implies X\text{ is quasi uniformly disconnected.}
\end{equation*}%
We will also compute the quasi-Assouad dimension for Moran set.

\begin{definition}
For any $\delta \in (0,1)$, let
\begin{equation*}
h_{X}(\delta )=\inf \{\alpha \geq 0:\exists b,c>0\text{ s.t. }N_{r,R}(X)\leq
c\left( \frac{R}{r}\right) ^{\alpha }\ \forall \ 0<r<r^{1-\delta }\leq R<b\}.
\end{equation*}%
Then the quasi-Assouad dimension $\dim _{qA}X$ is defined by
\begin{equation*}
\dim _{qA}X=\lim_{\delta \rightarrow 0}h_{X}(\delta ).
\end{equation*}
\end{definition}

It is easy to check that \newline
(a) $\dim _{qA}E\leq \dim _{qA}F$ if $E\subset F;$ \newline
(b) $\dim _{qA}(E\cup F)=\max (\dim _{qA}E,\dim _{qA}F)$; \newline
(c) $\dim _{qA}E=\dim _{qA}f(E)$ if $f$ is a bi-Lipschitz mapping.

Then we have the following proposition.

\begin{proposition}
Suppose quasi-Assouad dimension is defined as above. Then\newline
$(1)$ $\overline{\dim }_{B}X\leq \dim _{qA}X\leq \dim _{A}X;$\newline
$(2)$ $\dim _{qA}E=\dim _{qA}g(E)$ if $g$ is a quasi-Lipschitz mapping$;$%
\newline
$(3)$ $\dim _{qA}X=\lim_{\delta \rightarrow 0}\overline{\lim }_{r\rightarrow
0}\sup_{r^{1-\delta }\leq R<|X|}\frac{\log N_{r,R}(X)}{\log R-\log r};$%
\newline
$(4)$ If $\dim _{qA}X<1$, then $X$ is quasi uniformly disconnected.
\end{proposition}

\begin{remark}
There exists a Moran set $E$ satisfying $\overline{\dim }_{B}E<\dim
_{qA}E<\dim _{A}E$, see Example \ref{ex:di} in the next subsection.
\end{remark}

\subsection{Moran set}

$\ $

Some special cases of Moran sets were first studied by Moran \cite{Moran}.
The later works \cite{Hua, HuaLi, Wen1,Wen2} developed the theory on the
geometrical structures and dimensions of Moran sets systematically.

Suppose that $J$ is an initial closed interval of $\mathbb{R}.$ Let $%
\{n_{k}\}_{k\geq 1}$ be an integer sequence satisfying $n_{k}\geq 2$ for all
$k$. Suppose $c_{k}\in (0,1/n_{k}]$ for all $k.$ Denote $\mathcal{D}%
^{k}=\{i_{1}\cdots i_{k}:i_{t}\in \mathbb{N}\cap \lbrack 1,n_{t}]\ $for all $%
t\}$ and $\mathcal{D}^{0}=\{\emptyset \}$ with empty word $\emptyset .$ Let $%
J_{\emptyset }=J.$ Suppose for any $k\geq 1\ $and any $i_{1}\cdots
i_{k-1}\in \mathcal{D}^{k-1}$, $J_{i_{1}\cdots i_{k-1}1},\cdots
,J_{i_{1}\cdots i_{k-1}n_{k}}$ are closed subintervals of $J_{i_{1}\cdots
i_{k-1}}$ with their interiors pairwise disjoint, such that the ratio
\begin{equation}
\frac{|J_{i_{1}\cdots i_{k-1}j}|}{|J_{i_{1}\cdots i_{k-1}}|}=c_{k}\text{ for
all }1\leq j\leq n_{k}.  \label{bili}
\end{equation}%
Then we call the following compact set
\begin{equation}
F=\bigcap_{k=0}^{\infty }\bigcup_{i_{1}\cdots i_{k}\in \mathcal{D}%
^{k}}J_{i_{1}\cdots i_{k}}  \label{moran}
\end{equation}%
a \emph{Moran set} with the structure $(J,\{n_{k}\}_{k},\{c_{k}\}_{k}).$ We
denote $F\in \mathcal{M}(J,\{n_{k}\}_{k},\{c_{k}\}_{k}).$ By \cite{Wen2}, if
$\sup_{k}n_{k}<\infty ,$ then
\begin{equation}
\dim _{H}F=\underline{\lim }_{k\rightarrow \infty }\frac{\log n_{1}\cdots
n_{k}}{-\log c_{1}\cdots c_{k}}\text{, }\overline{\dim }_{B}F=\overline{\lim
}_{k\rightarrow \infty }\frac{\log n_{1}\cdots n_{k}}{-\log c_{1}\cdots c_{k}%
}  \label{wen222}
\end{equation}

\begin{remark}
The Moran structure is quite different from the self-similar structure. In
Moran structure, the relative positions\ of subintervals $\{J_{i_{1}\cdots
i_{k-1}j}\}_{j=1}^{n_{k}}$ in $J_{i_{1}\cdots i_{k-1}}$ can be variant,
$J_{i_{1}\cdots i_{k-1}i_{k}}$ and $J_{i_{1}\cdots i_{k-1}i_{k}^{\prime }}$
may share a common endpoint.
\end{remark}

If for any $k$ and any $i_{1}\cdots i_{k-1},J_{i_{1}\cdots i_{k-1}1},\cdots
,J_{i_{1}\cdots i_{k-1}n_{k}}$ are distributed uniformly in $J_{i_{1}\cdots
i_{k-1}}$ such that $J_{i_{1}\cdots i_{k-1}1}$ ($J_{i_{1}\cdots i_{k-1}n_{k}}$) and $J_{i_{1}\cdots i_{k-1}}$
share the left (right) endpoint, we call the Moran set a
uniform Cantor set.

For Moran set $F\in \mathcal{M}(J,\{n_{k}\},\{c_{k}\})$ with $%
\inf_{k}c_{k}>0 $, Li, Li, Miao and Xi \cite{LLMX} obtained the Assouad
dimension%
\begin{equation}
\dim _{A}F=\lim_{m\rightarrow \infty }\sup_{k}\frac{\log (n_{k+1}\cdots
n_{k+m})}{-\log (c_{k+1}\cdots c_{k+m})}.  \label{LLMX}
\end{equation}%
For uniform Cantor set $K$ with parameters $(J,\{n_{k}\}_{k},\{c_{k}\}_{k}),$
one recent result by Peng, Wang and Wen \cite{PWW}\ is that $\dim _{A}K=1$
if $\sup_{k}n_{k}=+\infty $.

\begin{theorem}
If $\lim_{k\rightarrow \infty }\frac{\log c_{k}}{\log c_{1}\cdots c_{k}}=0$,
then for any $E\in \mathcal{M}(J,\{n_{k}\},\{c_{k}\}),$
\begin{equation*}
\dim _{qA}E=\lim_{\delta \rightarrow 0}\underset{q\rightarrow \infty }{%
\overline{\lim }}\max_{1\leq p\leq l_{q,\delta }}\frac{\log (n_{p}\cdots
n_{q})}{-\log (c_{p}\cdots c_{q})},
\end{equation*}%
where $l_{q,\delta }=\max \{1\leq p\leq q:\frac{\log (c_{p}\cdots c_{q})}{%
\log (c_{1}\cdots c_{q})}>\delta \}$. In particular, if $\inf_{k\geq
1}c_{k}>0$, we obtain that
\begin{equation*}
\dim _{qA}E=\lim_{\eta \rightarrow 0}\underset{q\rightarrow \infty }{%
\overline{\lim }}\max_{1\leq p\leq q(1-\eta )}\frac{\log (n_{p}\cdots n_{q})%
}{-\log (c_{p}\cdots c_{q})}.
\end{equation*}
\end{theorem}

\begin{corollary}
If $\lim_{k\rightarrow \infty }\frac{\log c_{k}}{\log c_{1}\cdots c_{k}}=0$
and $\overline{\lim }_{k\rightarrow \infty }\frac{\log n_{k}}{-\log c_{k}}%
<1, $ then any $E\in \mathcal{M}(J,\{n_{k}\},\{c_{k}\})$ is quasi uniformly
disconnected.
\end{corollary}

The following Moran sets in Examples \ref{ex:inf}-\ref{ex:di2} are quasi
uniformly disconnected but not uniformly disconnected. In fact, for compact
set $E\subset \mathbb{R}^{1},$ $\dim _{A}E<1$ if and only if $E$ is
uniformly disconnected \cite{Luukkainen}.

\begin{example}
\label{ex:inf}Consider a uniform Cantor set $K$\ with $J=[0,1],$ $%
n_{k}=3^{k},c_{k}=3^{-2k}$ for all $k$. Then $\lim_{k\rightarrow \infty }%
\frac{\log c_{k}}{\log c_{1}\cdots c_{k}}=0$ and $\sup_{k}n_{k}=+\infty $,
using Theorem 1 and the result on Assouad dimension by Peng et al, we have
\begin{equation*}
\frac{1}{2}=\dim _{qA}K<\dim _{A}K=1.
\end{equation*}
\end{example}

\begin{example}
\label{ex:di}Let $\{q_{1}<q_{2}<\cdots <q_{t}<q_{t+1}<\cdots \}$ be a
positive integer sequence such that $q_{t+1}>2q_{t}$ for all $t$ and $%
\lim\limits_{t\rightarrow \infty }\frac{t}{q_{t}}=\lim\limits_{t\rightarrow
\infty }\frac{q_{1}+q_{2}+\cdots +q_{t-1}}{q_{t}}=0.$ Let $n_{k}\equiv 2$
and
\begin{equation*}
c_{k}=\left\{
\begin{array}{ll}
1/4 & \text{if }k\in (q_{t},2q_{t}], \\
(1-\frac{1}{2t})/2 & \text{if }k\in (2q_{t},2q_{t}+t], \\
1/5 & \text{otherwise.}%
\end{array}%
\right.
\end{equation*}%
According to (\ref{wen222}), Theorem 1 and (\ref{LLMX}), for any $E\in
\mathcal{M}([0,1],\{n_{k}\},\{c_{k}\})$ we have
\begin{equation*}
\overline{\dim }_{B}E=\frac{2\log 2}{\log 5+\log 4}<\dim _{qA}E=\frac{1}{2}%
<\dim _{A}E=1.
\end{equation*}
\end{example}

\begin{example}
\label{ex:di2}Let $\{q_{t}\}_{t}$ and $\{n_{k}\}_{k}$ be defined as in
Example \ref{ex:di}. Suppose $f:[1,2]\rightarrow (2,5)$ is a continuous
function$.$ Let
\begin{equation*}
\bar{c}_{k}=\left\{
\begin{array}{ll}
1/f(\frac{k}{q_{t}}) & \text{if }k\in (q_{t},2q_{t}], \\
(1-\frac{1}{2t})/2 & \text{if }k\in (2q_{t},2q_{t}+t], \\
1/5 & \text{otherwise.}%
\end{array}%
\right.
\end{equation*}%
Then for any $F\in \mathcal{M}([0,1],\{n_{k}\},\{\bar{c}_{k}\})$ we have
\begin{equation*}
\overline{\dim }_{B}F=\frac{\log 2}{(\log 5+\int\nolimits_{1}^{2}\log
f(x)dx)/2}<\dim _{qA}F=\frac{\log 2}{\log \left( \min\limits_{1\leq x\leq
2}f(x)\right) }<\dim _{A}F=1.
\end{equation*}
\end{example}

\subsection{Quasi-Lipschitz equivalence of Moran sets}

$\ $

We say $E\in \mathcal{M}(J,\{n_{k}\}_{k},\{c_{k}\}_{k})$ is of slow change,
if
\begin{equation}
\lim_{k\rightarrow \infty }\frac{\log c_{k}}{\log c_{1}\cdots c_{k}}=0\text{
and}\inf_{k}\frac{\log n_{k}}{\log c_{k}}>0.
\end{equation}%
In fact, if $\inf_{k}c_{k}>0,$ then $\inf_{k}\frac{\log n_{k}}{\log c_{k}}%
>0. $ The scale function $g_{E}(r)$ of $E$ is defined by
\begin{equation}
g_{E}(r)=\frac{\log n_{1}\cdots n_{k}}{-\log c_{1}\cdots c_{k}}\text{ if }%
c_{1}\cdots c_{k}\leq \frac{r}{|J|}<c_{1}\cdots c_{k-1}.
\end{equation}

It is proved in \cite{Lv2} that two quasi uniformly disconnected Moran sets $%
E,F$ of slow change are quasi-Lipschitz equivalent if and only if
\begin{equation}
\lim_{r\rightarrow 0}\frac{g_{E}(r)}{g_{F}(r)}=1.  \label{222}
\end{equation}

\begin{example}
Let $J=[0,1]$, $n_{k}\equiv 2,$ $c_{k}\equiv 1/4.$ Then by (\ref{222}), the
uniform Cantor set $K$ in Example \ref{ex:inf} and any $E\in \mathcal{M}%
([0,1],\{n_{k}\}_{k},\{c_{k}\}_{k})$ are quasi-Lipschitz equivalent,
although their structures seem to be quite different.
\end{example}

\begin{example}
Suppose $\{q_{t}\}_{t}$ and $\{n_{k}\}_{k}$ are given as in Example \ref%
{ex:di}. When $k\in (q_{t},2q_{t}+t]$ for some $t,$ we define $c_{k}=d_{k}$
as in Example \ref{ex:di}. When $k\in (2q_{t}+t,q_{t+1}]$ for some $t,$ we
can select $c_{k},d_{k}\ $from $\{1/5,1/6\}.$ Then $E\in \mathcal{M}%
([0,1],\{n_{k}\}_{k},\{c_{k}\}_{k}),F\in \mathcal{M}([0,1],\{n_{k}\}_{k},%
\{d_{k}\}_{k})$ are Moran sets of slow change, and they are quasi-Lipschitz
equivalent if and only if
\begin{equation*}
\lim_{k\rightarrow \infty }\frac{\#\{i\leq k:c_{i}=1/5\}}{\#\{i\leq
k:d_{i}=1/5\}}=1\text{ or }\lim_{k\rightarrow \infty }\frac{\#\{i\leq
k:c_{i}=1/6\}}{\#\{i\leq k:d_{i}=1/6\}}=1.
\end{equation*}
\end{example}

\medskip

The paper is organized as follows. We prove the basic properties of
quasi-Assouad dimension in Section 2. In particular, we use the idea \cite%
{MS}\ by Mattila and Saaranen to verify (4) of Proposition 1. In section 3,
we compute the quasi-Assouad dimension under the assumption $%
\lim_{k\rightarrow \infty }\frac{\log c_{k}}{\log c_{1}\cdots c_{k}}=0.$ In
the last section, we discuss the quasi-Assouad dimension for general Moran
sets.

\medskip

\section{Basic properties of quasi-Assouad dimension}

\subsection{Detail of Example 1}

$\ $

Note that $\sum_{i=1}^{\infty }(i+1)^{-\alpha }=\infty $, then $%
\lim_{k\rightarrow \infty }a_{k}=0$.

It is easy to check that $\{a_{i}-a_{i+1}\}_{i}$ is decreasing$.$ We will
use the following estimation \newline
$(\mathrm{i})$ $\lim_{k\rightarrow \infty }\frac{\log a_{k}}{\log a_{k-1}}%
=1, $ \newline
$(\mathrm{ii})$ $\lim_{k\rightarrow \infty }\frac{\log k}{\log a_{k}}=0,$
\newline
$(\mathrm{iii})$ $\lim_{k\rightarrow \infty }\frac{\log (a_{k}-a_{k+1})/2}{%
\log a_{k}}=\lim_{k\rightarrow \infty }\frac{\log (a_{k-1}-a_{k})/2}{\log
(a_{k}-a_{k+1})/2}=1.$

For $(\mathrm{i})$, $\lim_{k\rightarrow \infty }\frac{\log a_{k}}{\log
a_{k-1}}=1+\lim_{k\rightarrow \infty }\frac{\log (1-(k+1)^{-\alpha })}{\log
a_{k}}=1.$ For $(\mathrm{ii})$, by Stolz theorem and $(\mathrm{i})$, we have
\begin{eqnarray*}
\lim_{k\rightarrow \infty }\frac{\log k}{-\log a_{k}}=\lim_{k\rightarrow
\infty }\frac{\log (k+1)}{-\log a_{k-1}} &=&\lim_{k\rightarrow \infty }\frac{%
\log \frac{k+1}{k}}{-\log (a_{k-1}/a_{k-2})} \\
&=&\lim_{k\rightarrow \infty }\frac{\log (1+\frac{1}{k})}{-\log
(1-k^{-\alpha })}=\lim_{k\rightarrow \infty }k^{-1+\alpha }=0.
\end{eqnarray*}%
Since $a_{k}-a_{k+1}=(k+2)^{-\alpha }a_{k},$ using Stolz theorem and $(%
\mathrm{ii}),$ we obtain $(\mathrm{iii}).$

\medskip

(1) Firstly, the set $E$ is not uniformly disconnected and thus $\dim
_{A}E=1.$ Otherwise, suppose that there exist constant $0<c<1$ and $r^{\ast
}>0$ such that for any $x\in E$ and any $0<r<r^{\ast }$, we can find a set $%
E_{x,r}\subset E$ satisfying
\begin{equation*}
E\cap B(x,cr)\subset E_{x,r}\subset B(x,r)\text{~and~dist}%
(E_{x,r},E\backslash E_{x,r})\geq cr.
\end{equation*}%
Notice that the gap sequence $\{a_{i}-a_{i+1}\}_{i}$ is decreasing. Take $%
k\geq 1$ large enough such that $\frac{a_{k}}{2}<r^{\ast }$ and $\frac{%
a_{k}-a_{k+1}}{a_{k}}=(k+2)^{-\alpha }\leq \frac{c}{3},$ then for $x=a_{k}$
and $r=\frac{a_{k}}{2}$, we can not find such $E_{x,r}$. In fact, for any
subset $F$ of $E\cap B(x,r)$ containing $x$, since gap sequence is
decreasing, we must have
\begin{equation*}
\text{dist}(F,E\backslash F)\leq a_{k}-a_{k+1}\leq \frac{ca_{k}}{3}<cr.
\end{equation*}%
This is a contradiction.

\medskip

(2) Secondly, we will prove that $E$ is quasi uniformly disconnected.

For $0<r<1$, let $\psi (r)=\frac{a_{k}-a_{k+1}}{2}$ if $a_{k}\leq r<a_{k-1}$
for some $k\geq 1$ with $a_{0}=1$. Then $\lim_{r\rightarrow 0}\frac{\log
\psi (r)}{\log r}=1$ by $(\mathrm{iii}).$

Given any $x\in E$ and any $0<r<1$, suppose $a_{k}\leq r<a_{k-1}$ for some $%
k\geq 1$. If $0\leq x\leq a_{k+1}$, then we can take $E_{x,r}=[0,a_{k+1}]%
\cap E$; if $a_{k+1}<x\leq 1$, then we can take $E_{x,r}=\{x\}$. Then
\begin{equation*}
E\cap B(x,\psi (r))\subset E_{x,r}\subset B(x,r)\text{~and~dist}%
(E_{x,r},E\backslash E_{x,r})\geq \psi (r).
\end{equation*}

\medskip

(3) Furthermore, we can obtain that $\dim _{H}E=\dim _{B}E=\dim _{qA}E=0$.

It suffices to verify $h_{E}(\delta )\leq \varepsilon \ $for any fixed $%
\varepsilon >0$ and $\delta \in (0,1)$.

Given $a_{k}\leq r<a_{k-1}$ and $a_{l}\leq R<a_{l-1}$ with $r^{1-\delta
}\leq R,$ we consider $N_{r,R}.$ Suppose $m_{k}$ is an integer such that
\begin{equation*}
\frac{a_{m_{k}}-a_{m_{k}+1}}{2}\leq a_{k}<\frac{a_{m_{k}-1}-a_{m_{k}}}{2}.
\end{equation*}%
since $\{a_{i}-a_{i+1}\}_{i}$ is decreasing. Now we obtain
\begin{equation}
\lim_{k\rightarrow \infty }\frac{\log (a_{m_{k}}/a_{k})}{\log a_{k-1}}=0
\label{gggg1}
\end{equation}%
since $\lim\limits_{k\rightarrow \infty }\frac{\log a_{m_{k}}}{\log
(a_{m_{k}}-a_{m_{k}+1})/2}=\lim\limits_{k\rightarrow \infty }\frac{\log
(a_{m_{k}}-a_{m_{k}+1})/2}{\log a_{k-1}}=1$ and $\lim\limits_{k\rightarrow
\infty }\frac{\log a_{k}}{\log a_{k-1}}=1.$

By $(\mathrm{i})$-$(\mathrm{iii})$, we also have
\begin{equation}
\lim_{k\rightarrow \infty }\frac{\log 2m_{k}}{\log a_{k-1}}%
=\lim_{k\rightarrow \infty }\frac{\log m_{k}}{\log (a_{m_{k}}-a_{m_{k+1}})/2}%
=\lim_{m_{k}\rightarrow \infty }\frac{\log m_{k}}{\log a_{m_{k}}}=0.
\label{gggg2}
\end{equation}

Since
\begin{equation*}
E=([0,a_{m_{k}}]\cap E)\cup ([a_{m_{k}-1},1]\cap E),
\end{equation*}%
for $R$ small enough, using (\ref{gggg1}), (\ref{gggg2}) and $r^{1-\delta
}\leq R,$ we obtain that
\begin{equation*}
N_{r,R}(E)\leq N(E,a_{k})\leq \frac{a_{m_{k}}}{2a_{k}}+m_{k}\leq \max \{%
\frac{a_{m_{k}}}{a_{k}},2m_{k}\}\leq (a_{k-1})^{-\delta \varepsilon }\leq
r^{-\delta \varepsilon }\leq (\frac{R}{r})^{\varepsilon },
\end{equation*}%
where $N(E,a_{k})$ is the smallest number of balls with radii $a_{k}$ needed
to cover $E.$

\subsection{Proof of Proposition 1}

$\ $

When $X$ is fixed, we use $N_{r,R}$ and $N(r)$ to represent $N_{r,R}(X)$ and
$N_{r,|X|}(X)$ respectively.

\begin{proof}[Proof of (1) in Proposition 1]
$\ $

It is clear that $\dim _{qA}X\leq \dim _{A}X$. Now we shall verify that for
any $\delta \in (0,1)$,
\begin{equation*}
\overline{\dim }_{B}X\leq h_{X}(\delta ).
\end{equation*}%
In fact, for fixed $\delta \in (0,1)$, we can assume that for any $\alpha
>h_{X}(\delta )$ there are $b,c>0$ such that for $0<r<r^{1-\delta }\leq R<b$,%
\begin{equation*}
N_{r,R}\leq c\left( \frac{R}{r}\right) ^{\alpha }.
\end{equation*}%
Fix some $R<b$. When $r$ is small enough, using $N(r)\leq N_{r,R}\cdot N(R)$%
, we have%
\begin{equation*}
N(r)\leq N(R)\cdot c\left( \frac{R}{r}\right) ^{\alpha },
\end{equation*}%
which implies%
\begin{equation*}
\frac{\log N(r)}{-\log r}\leq \frac{\log N(R)}{-\log r}+\frac{\log c}{-\log r%
}+\frac{\alpha \log R}{-\log r}+\alpha .
\end{equation*}%
Letting $r\rightarrow 0,$ we obtain that $\overline{\dim }_{B}(X)\leq \alpha
$ and thus $\overline{\dim }_{B}(X)\leq h_{X}(\delta )$.
\end{proof}

\bigskip

\begin{proof}[Proof of (2) in Proposition 1]
$\ $

It suffices to show that $\dim _{qA}F\leq \dim _{qA}E.$

Suppose $g:E\rightarrow F$ is a quasi-Lipschitz mapping. Then there exist
increasing functions $\phi ,\zeta :\mathbb{R}^{+}\rightarrow \mathbb{R}^{+}$
with
\begin{equation}
\lim\limits_{r\rightarrow 0}\frac{\log \phi (r)}{\log r}=\lim\limits_{r%
\rightarrow 0}\frac{\log \zeta (r)}{\log r}=1  \label{555}
\end{equation}%
such that for all $y\in F,$%
\begin{equation}
B(y,R)\subset g(B(g^{-1}y,\phi (R)))\text{ and }g(B(g^{-1}y,\zeta
(r)))\subset B(y,r).  \label{666}
\end{equation}

Given $\delta \in (0,1/2)$ and any $\alpha _{\delta }>h_{E}(\delta ),$ for
all $r<r^{1-2\delta }\leq R$ with $R$ small enough, using (\ref{666}) we
obtain
\begin{equation*}
N_{r,R}(F)\leq N_{\zeta (r),\phi (R)}(E)\leq c_{\delta }(\frac{\phi (R)}{%
\zeta (r)})^{\alpha _{\delta }},
\end{equation*}%
where $c_{\delta }>0$ is a constant and $\zeta (r)^{1-\delta }\leq \phi (R)$
due to $r^{1-2\delta }\leq R$ and (\ref{555})$.$ Since $\frac{\log \phi (R)}{%
\log R}\rightarrow 1,$ $\frac{\log \zeta (r)}{\log r}\rightarrow 1$ and $0<%
\frac{\log R}{\log r}\leq 1-2\delta ,$ we have
\begin{equation*}
\frac{\log \phi (R)-\log \zeta (r)}{\log R-\log r}\rightarrow 1\text{ as }%
R\rightarrow 0,
\end{equation*}%
which implies that for any fixed $\varepsilon >0$,%
\begin{equation*}
N_{r,R}(F)\leq c_{\delta }(\frac{\phi (R)}{\zeta (r)})^{\alpha _{\delta
}}\leq c_{\delta }(\frac{R}{r})^{\alpha _{\delta }(1+\varepsilon )}
\end{equation*}%
whenever $R$ is small enough. Then means $h_{F}(2\delta )\leq h_{E}(\delta
). $ Letting $\delta \rightarrow 0,$ we have $\dim _{qA}F\leq \dim _{qA}E.$
\end{proof}

\bigskip

\begin{proof}[Proof of (3) in Proposition 1]
$\ $

We shall verify that for $\delta \in (0,1)$,
\begin{equation}
h_{X}(\delta )=\overline{\lim_{r\rightarrow 0}}\sup_{r^{1-\delta }\leq R<|X|}%
\frac{\log N_{r,R}}{\log R-\log r}.  \label{E:h}
\end{equation}

For any fixed number $\alpha >h_{X}(\delta )$, there exist $b$, $c>0$ such
that for all $0<r<r^{1-\delta }\leq R<b$,
\begin{equation*}
N_{r,R}\leq c\left( \frac{R}{r}\right) ^{\alpha }.
\end{equation*}%
When $r$ is small enough, if $R\geq b$, then $N_{r,R}\leq N_{r,\frac{b}{2}%
}\cdot N_{\frac{b}{2},R}\leq N_{r,\frac{b}{2}}\cdot N(\frac{b}{2}).$ That
means for $r^{1-\delta }\leq R<|X|,$ we obtain that%
\begin{equation*}
N_{r,R}\leq \left\{
\begin{array}{ll}
c\left( \frac{R}{r}\right) ^{\alpha } & \text{ if }R<b, \\
cN(\frac{b}{2})\left( \frac{b}{2r}\right) ^{\alpha } & \text{ if }R\geq b.%
\end{array}%
\right.
\end{equation*}%
Therefore,
\begin{align*}
& \sup_{r^{1-\delta }\leq R<|X|}\frac{\log N_{r,R}}{\log R-\log r} \\
\leq & \alpha +\max \{\sup_{r^{1-\delta }\leq R<b}\frac{\log c}{\log R-\log r%
},\sup_{b\leq R}\frac{\log c+\log N(\frac{b}{2})-\alpha \log 2}{\log b-\log r%
}\},
\end{align*}%
which implies $\overline{\lim }_{r\rightarrow 0}\sup_{r^{1-\delta }\leq
R<|X|}\frac{\log N_{r,R}}{\log R-\log r}\leq \alpha $ for any $\alpha
>h_{X}(\delta )$. Hence
\begin{equation*}
\overline{\lim_{r\rightarrow 0}}\sup_{r^{1-\delta }\leq R<|X|}\frac{\log
N_{r,R}}{\log R-\log r}\leq h_{X}(\delta ).
\end{equation*}

On the other hand, for any $\alpha _{0}>\overline{\lim }_{r\rightarrow
0}\sup_{r^{1-\delta }\leq R<|X|}\frac{\log N_{r,R}}{\log R-\log r},$ there
exists $r_{0}\in (0,1)$ such that for $r<r_{0}$,%
\begin{equation*}
\alpha _{0}>\sup_{r^{1-\delta }\leq R<|X|}\frac{\log N_{r,R}}{\log R-\log r},
\end{equation*}%
Hence for $0<r<r^{1-\delta }\leq R<r_{0}$, we obtain that
\begin{equation*}
\frac{\log N_{r,R}}{\log R-\log r}<\alpha _{0},
\end{equation*}%
which implies $N_{r,R}\leq \left( \frac{R}{r}\right) ^{\alpha _{0}}.$
Therefore $\alpha _{0}\geq h_{X}(\delta )$, and thus
\begin{equation*}
\overline{\lim_{r\rightarrow 0}}\sup_{r^{1-\delta }\leq R<|X|}\frac{\log
N_{r,R}}{\log R-\log r}\geq h_{X}(\delta ).
\end{equation*}
\end{proof}

\bigskip

We will use the idea of \cite{MS} by Mattila and Saaranen.

\begin{proof}[Proof of (4) in Proposition 1]
\

Suppose $\dim _{qA}X=s<1.$ Without loss of generality, we assume that
\begin{equation*}
|X|<1.
\end{equation*}%
For any $\delta \in (0,1)$ and $r\in (0,1)$, let
\begin{equation*}
h_{X}(\delta ,r)=\sup_{r^{1-\delta }\leq R<|X|}\frac{\log N_{r,R}}{\log
R-\log r},
\end{equation*}%
Then by (\ref{E:h}), we have $\overline{\lim }_{r\rightarrow 0}h_{X}(\delta
,r)=h_{X}(\delta )$.

For any $\delta \in (0,1)$, there exists $r_{\ast }(\delta )\in (0,1)$ such
that for all $r\in (0,r_{\ast }(\delta ))$,
\begin{equation*}
h_{X}(\delta ,r)<h_{X}(\delta )+\frac{1-s}{2}\leq s+\frac{1-s}{2}=\frac{1}{2}%
(1+s).
\end{equation*}%
Hence when $r\in (0,r_{\ast }(\delta ))$ and $r^{1-\delta }\leq R<|X|$, we
always have
\begin{equation}
N_{r,R}\leq \left( \frac{R}{r}\right) ^{\frac{1}{2}(1+s)}.  \label{E:lem3N}
\end{equation}

For any $0<R<|X|$, let
\begin{equation*}
\rho (R)=\inf_{0<\delta <1}\max \{\frac{\log r_{\ast }(\delta )}{\log R}-1,%
\frac{\delta }{\delta -1},\frac{2\log 6}{-(1-s)\log R}\}.
\end{equation*}%
Then we have
\begin{equation*}
\rho (R)\geq \frac{2\log 6}{-(1-s)\log R}>0.
\end{equation*}%
Fix arbitrary $\delta \in (0,1)$,
\begin{equation*}
\underset{R\rightarrow 0}{\overline{\lim }}\rho (R)\leq \underset{%
R\rightarrow 0}{\overline{\lim }}\max \{\frac{\log r_{\ast }(\delta )}{\log R%
}-1,\frac{\delta }{\delta -1},\frac{2\log 6}{-(1-s)\log R}\}\leq \frac{%
\delta }{\delta -1}.
\end{equation*}%
Letting $\delta \rightarrow 0,$ we have
\begin{equation*}
\lim_{R\rightarrow 0}\rho (R)=0.
\end{equation*}

Fix $x\in X$ and $R$ small enough, let $B_{0}=\{x\}$ and
\begin{equation*}
B_{i}=B(x,iR^{1+2\rho (R)})\backslash B(x,(i-1)R^{1+2\rho (R)})\text{ for
any }i\in \mathbb{N}\text{.}
\end{equation*}%
Write $n_{R}=[R^{-2\rho (R)}]$, where $[z]$ denotes the integral part of $z$%
. Then $n_{R}\geq \lbrack R^{\frac{4\log 6}{(1-s)\log R}}]=[6^{\frac{4}{1-s}%
}]>1$.

We conclude that

\begin{claim}
If $R$ is small enough, then there exists $1\leq i\leq n_{R}$ such that
\begin{equation*}
B_{i}\cap X=\varnothing .
\end{equation*}
\end{claim}

Otherwise, we assume that for all $1\leq i\leq n_{R}$, $B_{i}\cap X$ is
non-empty. We take $y_{i}\in B_{i}\cap X$ and let
\begin{equation*}
\Theta =\{y_{i}:1\leq i\leq n_{R}\}.
\end{equation*}

Write $r=R^{1+2\rho (R)}.$ For any $1\leq i,j\leq n_{R}$ with $j\geq i+3$,
the closed ball
\begin{eqnarray*}
B(y_{j},r) &\subset &X\backslash B(x,(j-2)r), \\
B(y_{i},r) &\subset &B(x,(i+1)r),
\end{eqnarray*}%
which implies $B(y_{j},r)\cap B(y_{i},r)=\varnothing ,$ then any ball with
radius $R^{1+2\rho (R)}$ will cover at most $3$ points in $\Theta $, i.e.,
\begin{equation}
\frac{n_{R}}{3}\leq N_{R^{1+2\rho (R)},R}.  \label{E:nR3}
\end{equation}

By the definition of $\rho (R)$, there exists $0<\delta <1$ such that%
\begin{equation*}
\max \{\frac{\log r_{\ast }(\delta )}{\log R}-1,\frac{\delta }{\delta -1},%
\frac{2\log 6}{-(1-s)\log R}\}<2\rho (R),
\end{equation*}%
Therefore
\begin{equation*}
R^{1+2\rho (R)}<r_{\ast }(\delta )\text{, }R^{(1+2\rho (R))}<R^{(1+2\rho
(R))(1-\delta )}<R,
\end{equation*}%
and
\begin{equation}
\frac{\log 6}{-(1-s)\log R}<\rho (R).  \label{ai}
\end{equation}

Let $r=R^{1+2\rho (R)}$, by (\ref{E:lem3N}) we have
\begin{equation*}
N_{R^{1+2\rho (R)},R}\leq R^{-(1+s)\rho (R)}.
\end{equation*}%
Using (\ref{E:nR3}) and the definition of $n_{R}$, we obtain that%
\begin{equation*}
\frac{R^{-2\rho (R)}}{6}\leq \frac{n_{R}}{3}\leq R^{-(1+s)\rho (R)},
\end{equation*}%
which implies
\begin{equation*}
\rho (R)\leq \frac{\log 6}{-(1-s)\log R}.
\end{equation*}%
This is contradictory to (\ref{ai}). The claim is proved.

\medskip

Using Claim 1, take $\psi (R)=R^{1+2\rho (R)}$ and $1\leq i\leq n_{R}$ such
that $B_{i}\cap X$ is non-empty, then let
\begin{equation*}
E_{x,R}=X\cap B(x,(i-1)\psi (R)),
\end{equation*}%
satisfying
\begin{equation*}
X\cap B(x,\psi (R))\subset E_{x,R}\subset B(x,R)\text{ and d}%
_{X}(E_{x,R},X\backslash E_{x,R})\geq \psi (R).
\end{equation*}%
Therefore the quasi uniform disconnectedness is obtained.
\end{proof}

\bigskip

\section{Quasi-Assouad dimension of Moran set}

We call the closed interval $J_{i_{1}\cdots i_{k}}$ a basic interval of rank
$k.$ For $1\leq p\leq q$, let
\begin{equation}
s_{p,q}=\frac{\log (n_{p}\cdots n_{q})}{-\log (c_{p}\cdots c_{q})}.
\label{pq}
\end{equation}

\begin{proof}[Proof of Theorem 1]
$\ $

Notice that $\dim _{qA}E=\lim_{\delta \rightarrow 0}h_{E}(\delta )$. It
suffices to show
\begin{equation*}
h_{E}(\delta )=\overline{\lim }_{q\rightarrow \infty }\max_{1\leq p\leq
l_{q,\delta }}s_{p,q},\text{ for any }\delta \in (0,1).
\end{equation*}

\textbf{Step 1.} We shall verify $h_{E}(\delta )\leq \overline{\lim }%
_{q\rightarrow \infty }\max_{1\leq p\leq l_{q,\delta }}s_{p,q},.$

For any $s>\overline{\lim }_{q\rightarrow \infty }\max_{1\leq p\leq
l_{q,\delta }}s_{p,q}$, we will find $b,c>0$ such that for all $%
0<r<r^{1-\delta }\leq R<b,$
\begin{equation*}
N_{r,R}\leq c\left( \frac{R}{r}\right) ^{s}.
\end{equation*}

Since $s>\overline{\lim }_{q\rightarrow \infty }\max_{1\leq p\leq
l_{q,\delta }}s_{p,q},$ we can find a constant $\sigma >0$ small enough such
that $s(1-\sigma )>\overline{\lim }_{q\rightarrow \infty }\max_{1\leq p\leq
l_{q,\delta }}s_{p,q}$. Then there exists $N\in \mathbb{N}$ such that for
all $q\geq N$ and all $k\geq N,$%
\begin{equation}
s(1-\sigma )>s_{p,q}\text{ for any }1\leq p\leq l_{q,\delta },  \label{wwww}
\end{equation}%
and
\begin{equation*}
\frac{\log c_{k}}{\log c_{1}\cdots c_{k}}<\frac{\delta \sigma }{2}.
\end{equation*}

Take $b=c_{1}\cdots c_{N}$ and $c=3$. For all $0<r<r^{1-\delta }\leq R<b$,
we assume that
\begin{equation*}
c_{1}\cdots c_{p}\leq R<c_{1}\cdots c_{p-1}\text{~and~}c_{1}\cdots c_{q}\leq
r<c_{1}\cdots c_{q-1},
\end{equation*}%
where $q\geq p\geq N+1$. Any ball of radius $R$ will intersect at most $3$
basic interval of rank $(p-1),$ i.e.,
\begin{equation*}
N_{r,R}\leq 3n_{p}\cdots n_{q}\leq 3\left( \frac{R}{r}\right) ^{\frac{\log
n_{p}\cdots n_{q}}{-\log c_{p}\cdots c_{q}+\log c_{p}c_{q}}}.
\end{equation*}

Since $(c_{1}\cdots c_{q})^{1-\delta }\leq r^{1-\delta }\leq R<c_{1}\cdots
c_{p-1}$, we have%
\begin{equation*}
\frac{\log c_{p}\cdots c_{q}}{\log c_{1}\cdots c_{q}}>\delta .
\end{equation*}%
That means $1\leq p\leq l_{q,\delta }$, using (\ref{wwww}) we have%
\begin{equation*}
s_{p,q}=\frac{\log n_{p}\cdots n_{q}}{-\log c_{p}\cdots c_{q}}<s(1-\sigma ).
\end{equation*}

On the other hand,
\begin{eqnarray*}
\frac{\log c_{p}c_{q}}{\log c_{p}\cdots c_{q}} &=&\frac{\log c_{p}+\log c_{q}%
}{\log c_{1}\cdots c_{q}}\cdot \frac{\log c_{1}\cdots c_{q}}{\log
c_{p}\cdots c_{q}} \\
&\leq &(\frac{\log c_{p}}{\log c_{1}\cdots c_{p}}+\frac{\log c_{q}}{\log
c_{1}\cdots c_{q}})\cdot \frac{1}{\delta } \\
&<&2\cdot \frac{\delta \sigma }{2}\cdot \frac{1}{\delta }=\sigma .
\end{eqnarray*}%
We obtain
\begin{equation*}
\frac{\log n_{p}\cdots n_{q}}{-\log c_{p}\cdots c_{q}+\log c_{p}c_{q}}\leq
\frac{\log n_{p}\cdots n_{q}}{-(1-\sigma )\log c_{p}\cdots c_{q}}=\frac{%
s_{p,q}}{1-\sigma }<s.
\end{equation*}%
Therefore,
\begin{equation*}
N_{r,R}\leq 3\left( \frac{R}{r}\right) ^{s}=c\left( \frac{R}{r}\right) ^{s}.
\end{equation*}

\textbf{Step 2.} We shall verify $h_{E}(\delta )\geq \overline{\lim }%
_{q\rightarrow \infty }\max_{1\leq p\leq l_{q,\delta }}s_{p,q}$ for any $%
\delta \in (0,1)$.

Fix $\delta \in (0,1)$, we assume that $\alpha >h_{E}(\delta )$, then there
exist $b,c>0$ such that for all $0<r<r^{1-\delta }\leq R<b$,
\begin{equation*}
N_{r,R}\leq c\left( \frac{R}{r}\right) ^{\alpha }.
\end{equation*}%
Without loss of generality, we assume that there exists $M\in \mathbb{N}$
such that
\begin{equation*}
c_{1}\cdots c_{M}<b\leq c_{1}\cdots c_{M-1}.
\end{equation*}

For any fixed number $\varepsilon \in (0,\frac{2(1-\delta )}{\delta })$, we
can take $N\in \mathbb{N}$ large enough such that for all $q\geq N\geq M$,
\begin{equation*}
\frac{\log c_{1}\cdots c_{M}}{\log c_{1}\cdots c_{q}}<\frac{\varepsilon }{%
2+\varepsilon }\text{~and~}\frac{\log 4c}{-\log c_{1}\cdots c_{q}}<\frac{%
\delta \varepsilon }{2}.
\end{equation*}%
For any $q\geq N$, let $r=c_{1}\cdots c_{q}$. Since
\begin{equation*}
\frac{\log c_{M+1}\cdots c_{q}}{\log c_{1}\cdots c_{q}}=1-\frac{\log
c_{1}\cdots c_{M}}{\log c_{1}\cdots c_{q}}>1-\frac{\varepsilon }{%
2+\varepsilon }>\delta ,
\end{equation*}%
then $M+1\leq l_{q,\delta }.$ For any integer $p\in \lbrack M+1,l_{q,\delta
}]$, let $R=c_{1}\cdots c_{p-1}$. Then $0<r<r^{1-\delta }\leq R<b$, and thus
\begin{equation*}
\frac{n_{p}\cdots n_{q}}{4}\leq N_{r,R}\leq c\left( \frac{R}{r}\right)
^{\alpha }=c\left( \frac{1}{c_{p}\cdots c_{q}}\right) ^{\alpha },
\end{equation*}%
which implies for $M<p\leq l_{j,q},$
\begin{align*}
\frac{\log n_{p}\cdots n_{q}}{-\log c_{p}\cdots c_{q}}& \leq \frac{\log 4c}{%
-\log c_{p}\cdots c_{q}}+\alpha \\
& =\frac{\log 4c}{-\log c_{1}\cdots c_{q}}\cdot \frac{\log c_{1}\cdots c_{q}%
}{\log c_{p}\cdots c_{q}}+\alpha <\frac{\delta \varepsilon }{2}\cdot \frac{1%
}{\delta }+\alpha =\alpha +\frac{\varepsilon }{2}.
\end{align*}

When $p\leq M$, using $n_{k}\leq c_{k}^{-1}$ for all $k,$ we obtain
\begin{align*}
\frac{\log n_{p}\cdots n_{q}}{-\log c_{p}\cdots c_{q}}& =\frac{\log
n_{p}\cdots n_{M}}{-\log c_{p}\cdots c_{q}}+\frac{\log n_{M+1}\cdots n_{q}}{%
-\log c_{p}\cdots c_{q}} \\
& \leq \frac{\log c_{p}\cdots c_{M}}{\log c_{p}\cdots c_{q}}+(\alpha +\frac{%
\varepsilon }{2}) \\
& \leq \frac{\log c_{1}\cdots c_{M}}{\log c_{M+1}\cdots c_{q}}+(\alpha +%
\frac{\varepsilon }{2}) \\
& =\frac{\log c_{1}\cdots c_{M}}{\log c_{1}\cdots c_{q}}\cdot \frac{\log
c_{1}\cdots c_{q}}{\log c_{M+1}\cdots c_{q}}+(\alpha +\frac{\varepsilon }{2})
\\
& <\frac{\varepsilon }{2+\varepsilon }\cdot \frac{2+\varepsilon }{2}+(\alpha
+\frac{\varepsilon }{2})=\alpha +\varepsilon .
\end{align*}

Therefore
\begin{equation*}
\max_{1\leq p\leq l_{q,\delta }}s_{p,q}\leq \alpha +\varepsilon .
\end{equation*}%
Letting $\varepsilon \rightarrow 0$, we obtain
\begin{equation*}
\max_{1\leq p\leq l_{q,\delta }}s_{p,q}\leq \alpha ,
\end{equation*}%
and thus $\overline{\lim }_{q\rightarrow \infty }\max_{1\leq p\leq
l_{q,\delta }}s_{p,q}\leq \alpha $. Therefore
\begin{equation*}
\overline{\lim }_{q\rightarrow \infty }\max_{1\leq p\leq l_{q,\delta
}}s_{p,q}\leq h_{E}(\delta )
\end{equation*}
follows$.$
\end{proof}

\section{Result of general Moran set}

Consider general Moran set $E\in \mathcal{M}\{J,\{n_{k}\}_{k},\{c_{k,i}\}_{k%
\geq 1,1\leq i\leq n_{k}}\}$ in Euclidean space $\mathbb{R}^{d}$ with $d\geq
1$, see \cite{Wen2} for details of general Moran set.

For any $p\leq q,$ let $s_{p,q}$ be the unique positive solution of the
equation $\Delta _{p,q}(s)=1$, where
\begin{equation}
\Delta _{p,q}(s)=\prod\nolimits_{i=p}^{q}\left(
\sum\nolimits_{j=1}^{n_{i}}(c_{i,j})^{s}\right) .  \label{eqns}
\end{equation}%
In fact, if $c_{k,1}=\cdots =c_{k,n_{k}}=c_{k}$ for all $k,$ then $s_{p,q}=%
\frac{\log n_{p}\cdots n_{q}}{-\log c_{p}\cdots c_{q}}$ as in (\ref{pq}).
Using the method in \cite{LLMX}, we can obtain the following result and skip
its proof.

\begin{proposition}
\label{P:general}If $\inf_{k,i}c_{k,i}>0$, then for any $E\in \mathcal{M}%
\{J,\{n_{k}\}_{k},\{c_{k,i}\}_{k\geq 1,1\leq i\leq n_{k}}\},$ we obtain that
\begin{equation*}
\dim _{qA}E=\lim_{\eta \rightarrow 0}\underset{q\rightarrow \infty }{%
\overline{\lim }}\max_{1\leq p\leq q(1-\eta )}s_{p,q}.
\end{equation*}
\end{proposition}

When $c_{k,1}=\cdots =c_{k,n_{k}}=c_{k}$ for all $k$ and $\inf_{k}c_{k}>0,$
the result is the same as in Theorem 1 since $s_{p,q}=\frac{\log n_{p}\cdots
n_{q}}{-\log c_{p}\cdots c_{q}}.$


\bigskip

\begin{thebibliography}{99}
\bibitem{Assouad1} P. Assouad, Espaces m\'{e}triques, plongements, facteurs,
Th\`{e}se de doctorat, Publ. Math. Orsay No. 223-7769, Univ. Paris XI,
Orsay, 1977.

\bibitem{Assouad2} P. Assouad, \`{E}tude d'une dimension m\'{e}trique li\'{e}%
e \`{a} la possibilit\'{e} de plongements dans $\mathbf{R}^{n}$, C.~R. Acad.
Sci. Paris S\'{e}r. A-B, 288(15), 1979, 731--734.

\bibitem{Assouad3} P. Assouad, Pseudodistances, facteurs et dimension m\'{e}%
trique, in Seminaire D'Analyse Harmonique (1979-1980), pp. 1--33, Publ.
Math. Orsay 80, 7, Univ. Paris XI, Orsay, 1980.

\bibitem{David Semmes} G. David, S. Semmes, Fractured fractals and broken
dreams: self-similar geometriy through metric and measure, Oxford University
Press, New York, 1997.

\bibitem{Edgar} G. Edgar, Integral, Probability, and Fractal Measures,
Springer-Verlag, New York, 1998.

\bibitem{Falconer} K J. Falconer, Fractal Geometry. Mathematical Foundations
and Applications, John Wiley Sons, Ltd., Chichester, 1990.

\bibitem{Fraser} J. M. Fraser, Assouad type dimensions and homogeneity of
fractals, arXiv:1301.2934, 2013.

\bibitem{Heinonen} J. Heinonen, Lectures on Analysis on Metric Spaces,
Springer-Verlag, New York, 2001.

\bibitem{Hua} S. Hua, On the Hausdorff dimension of generalized self-similar
sets (Chinese), Acta Math. Appl. Sinica, 17(4), 1994, 551--558.

\bibitem{HuaLi} S. Hua, W. X. Li, Packing dimension of generalized Moran
Sets, Progr. Natur. Sci. (English Ed.), 6(2), 1996, 148--152.

\bibitem{Lehrb} J. Lehrb, H. Tuominen, A note on the dimensions of Assouad
and Aikawa, J. Math. Soc. Japan, 65(2), 2013, 343--356.

\bibitem{Li} J. J. Li, Assouad dimensions of Moran sets, C. R. Math. Acad.
Sci. Paris, 351(1-2), 2013, 19--22.

\bibitem{LLMX} W. W. Li, W. X. Li, J. J. Miao, L. F. Xi, Assouad dimensions
of Moran sets and Cantor-like sets, ArXiv: 1404.4409, 2014.

\bibitem{Luukkainen} J. Luukkainen, Assouad dimension: Antifractal
metrization, porous sets, and homogeneous measures, J. Korean Math. Soc.,
35, 1998, 23--76.

\bibitem{Lv} F. L\"{u}, M. L. Lou, Z. Y. Wen, L. F. Xi, Bilipschitz
embedding of homogeneous set, arXiv:1402.0080, 2014.

\bibitem{Lv2} F. L\"{u}, M. L. Lou, L. F. Xi, Homogeneous sets and their
quasi-Lipschitz equivalence, Sci. China Math., in press, 2014.

\bibitem{Mackay} J. M. Mackay, Assouad dimension of self-affine carpets,
Conform. Geom. Dyn., 15, 2011, 177--187.

\bibitem{MT} J. M. Mackay, J. T. Tyson, Conformal dimension: Theory and
application. University Lecture Series, 54. American Mathematical Society,
Providence, RI, 2010.

\bibitem{Mattila} P. Mattila, Geometry of Sets and Measure in Euclidean
Spaces, Cambridge University Press, Cambridge, 1995.

\bibitem{MS} P. Mattila, P. Saaranen, Ahlfors--David regular sets and
bilipschitz maps, Ann. Acad. Sci. Fenn. Math., 34, 2009, 487--502.

\bibitem{Mcm} C. McMullen, The Hausdorff dimension of general Sierpinski
carpets, Nagoya Math. J. 96, 1-9, 1984.

\bibitem{Moran} P. A. Moran, Additive functions of intervals and Hausdorff
measure, Proc. Camb. Phil. Soc., 42, 1946, 15--23.

\bibitem{Olsen} L. Olsen, On the Assouad dimension of graph directed Moran
fractals, Fractals, 19, 2011, 221--226.

\bibitem{PWW} F. J. Peng, W. Wang and S. Y. Wen, On Assouad dimension of
products, preprint.

\bibitem{Tricot} C. Tricot, Curves and Fractal Dimension, Springer-Verlag,
New York, 1995.

\bibitem{Wang Xi N} Wang, Q. and Xi, L. F.: Quasi-Lipschitz equivalence of
Ahlfors-David regular sets, Nonlinearity, 24(3), 2011, 941-950.

\bibitem{Wang Xi S} Wang, Q. and Xi, L. F.: Quasi-Lipschitz equivalence of
quasi Ahlfors-David regular sets, Sci. China Math., 54(12), 2011, 2573--2582.

\bibitem{Wen1} Z. Y. Wen, Mathematical Foundations of Fractal Geometry,
Shanghai Scientific and Tech- nological Education Publishing House,
Shanghai, 2000.

\bibitem{Wen2} Z. Y. Wen, Moran sets and Moran classes, Chinese Sci. Bull.,
46, 2001, 1849--1856.\qquad
\end{thebibliography}
\end{document}